\documentclass[12pt]{article} \usepackage{amsmath,amsthm, amssymb} 
\usepackage{amssymb,latexsym}

\newtheorem*{theorem}{Theorem}

\begin{document}

\baselineskip=17pt 

\title{\bf On the remainder term in the circle problem in an arithmetic progression}

\author{\bf D. I. Tolev\footnote{Supported by the Ministry of Science and Education of Bulgaria, 
Grant DD VU 02/90.}}

\date{}
\maketitle

\centerline
{\bf \large Dedicated to 75th birthday of Professor A. A. Karatsuba }

\bigskip

\begin{abstract}
In this paper we improve the estimate for the remainder term 
in the asymptotic formula concerning the circle problem in an arithmetic progression.

\bigskip

Keywords: Circle problem, Arithmetic progressions.

Mathematics Subject Classification (2010): 11L05, 11N37.

\end{abstract}

\section{Introduction and statement of the result.}

\indent

Consider the sum
\begin{equation} \label{10}
  S_{q, a}(x) = \sum_{\substack{n \le x \\ n \equiv a \pmod {q} }} r(n) ,
\end{equation}
where $r(n)$ is the number of representations of $n$ as a sum of two squares.
In 1968 Smith~\cite{Smith} established that if $q = O \left( x^{\frac{2}{3}} \right) $ 
then we have the following asymptotic formula
\begin{equation} \label{12}
  S_{q, a} (x) = \pi  \frac{\eta_a(q)}{q^2} x  + 
  R_{q, a} (x) .
\end{equation}
The quantity $\eta_a(q) $ in the main term is defined by
\begin{equation} \label{17}
   \eta_a(q) = \# \{ \; 1 \le \alpha, \beta \le q  \; : \; \alpha^2 + \beta^2 \equiv a \pmod{q} \; \} 
\end{equation}
(the main term in Smith's paper is written in a slightly different form)
and $R_{q, a}(x)$ is the remainder term 
for which it is established that
\begin{equation} \label{17.5}
   R_{q, a} (x) \ll  x^{\frac{2}{3} + \xi} q^{-\frac{1}{2} (1 + 3 \xi) } \; 
   (q, a)^{\frac{1}{2}} \; \tau(q)  .
\end{equation}
Here $\tau(q)$ is the divisor function, $(q, a)$ stands for the greatest common divisor of $a$ and $q$
and $ 0 < \xi < \frac{1}{3}$.

\bigskip

We note that uniformly for $a$ we have
\begin{equation} \label{18}
   \eta_a(q) \ll q \, \tau(q) 
\end{equation}
(a proof is available in \cite[Lemma 2.8]{BloBruD})
and obviously $S_{q, a}(x) \ll x^{1+ \varepsilon} q^{-1}$ for any $\varepsilon > 0$.
Hence the asymptotic formula \eqref{12} 
is non-trivial for $q \le x^{\frac{2}{3} - \varepsilon}$.
However if $q$ is small then the estimate \eqref{17.5} for the error term is quite weak. 

\bigskip

In 1970 Varbanets~\cite{Varbanets}\footnote{The author would like to thank Mr Ping Xi for 
informing him about this paper.} 
considered the case $a=1$ and 
found better estimates for the remainder term in \eqref{12}.
He established that if  $x^{\frac{1}{2}} \le q \le x^{\frac{2}{3}}$ then
\begin{equation} \label{18.5}
  R_{q, 1}(x)  \ll \left( q^{\frac{1}{2} } \; + \;   x^{\frac{1}{2}} \, q^{- \frac{1}{4}} \right) \, x^{\varepsilon} .
\end{equation}
Varbanets also proved that if $ q \le x^{\frac{2}{3}} $ and 
\begin{equation} \label{18.7}
  q^{\frac{1}{4}} > \xi(q) , \qquad \text{where} \qquad \xi(q) = \prod_{p \mid q} p 
\end{equation}
(the product is taken over the prime divisors of $q$), then 
\begin{equation} \label{19}
  R_{q, 1}(x)  \ll \left( q^{\frac{1}{2} } \; + \;   x^{\frac{1}{2}} \, q^{- \frac{1}{2}} \, \xi(q) \right) \, x^{\varepsilon} .
\end{equation}
(The results of \cite{Varbanets} are actually slightly stronger --- with a power of $\tau(q)$ rather than
$x^{\varepsilon}$).

\bigskip

It is clear that our sum $S_{q, a}(x)$ is very similar to
\[
  T_{q, a}(x) = \sum_{\substack{n \le x \\ n \equiv a \pmod {q} }} \tau(n) .
\]
An asymptotic formula 
for $ T_{q, a}(x) $,
which is non-trivial for $q \le x^{\frac{2}{3} - \varepsilon}$, was established independently by 
A. Selberg and C. Hooley (in unpublished manuscripts). It is well known today that
if $(q, a)=1$ then we have
\begin{equation} \label{20}
  T_{q, a}(x) = \frac{\varphi(q)}{q^2} \, x 
  \left( 
   \log x + 2 \gamma - 1 + 2 \sum_{p \mid q} \frac{\log p}{p-1} \right)
   \,   + O \left( \tau^2(q) \, \big( q^{\frac{1}{2}} + x^{\frac{1}{3}} \big) \, \log x \right) ,
\end{equation}
where $\varphi(q)$ is the Euler function, $\gamma$ stands for the Euler constant
and where the summation is taken over the prime divisors of $q$.
A proof of \eqref{20} can be found for example in a recent paper of Blomer~\cite{Blomer} 
or in the book of Iwaniec and Kowalski \cite[Chapter 4]{IwKo}.
One may expect that a formula similar to \eqref{20} holds also for the sum $S_{q,a}(x)$.

\bigskip

In the present paper we prove the following:
\begin{theorem}
For the quantity $R_{q, a}(x)$ defined by \eqref{12} we have
\begin{equation} \label{30}
  R_{q, a}(x) \ll 
   \big( q^{\frac{1}{2} } + x^{ \frac{1}{3} } \big) \, (a, q)^{\frac{1}{2}} \, \tau^4(q) \, \log^4 x   .
\end{equation}
\end{theorem}

\bigskip

The asymptotic formula \eqref{12} with the estimate \eqref{30} for the remainder term
is non-trivial for $q \le x^{\frac{2}{3} - \varepsilon}$ and \eqref{30}
is stronger than \eqref{18.5} for all such $q$.
If we compare our 
bound for $R_{q, a}(x)$
with \eqref{19}, which is established only for $a=1$ and for $q$ satisfying \eqref{18.7}, 
we can see that \eqref{30} 
is stronger provided that $q \ll x^{\frac{1}{3}} \, \xi(q)^2$.

\bigskip

One of the main points in our proof is the estimation of the sum
$\mathcal H_{h, n}(q, a)$ defined by \eqref{340}. 
In Section~\ref{S1_b} we represent it as a linear combination of 
Kloosterman sums and then apply A. Weil's bound. 
An estimate of the same strength for $\mathcal H$ in the case $a=1$ 
is established in \cite{Varbanets} but appealing to a result of Bombieri~\cite{Bom}.
It is clear that Varbanets' method can be applied also for
the estimation of $\mathcal H$ for any integer $a$, but here we present
our method, which may be of some use in other occasions. 
The other arguments of the proof are elementary or based on the simplest 
theorems from the theory of the exponential sums.

\bigskip

Working in the same manner one may find another proof of \eqref{20}
and one may also establish similar asymptotic formulas for the quantities
\[
  \sum_{\substack{m^2 + n^2 \le x  \\ mn \equiv a \; (q) }} 1 \;
  \qquad \text{and} \qquad
  \sum_{\substack{mn \le x  \\ m^2 + n^2 \equiv a \; (q) } } 1 
\]
(here and later we write for simplicity $k \equiv a \; (q)$ instead of $k \equiv a \pmod{q}$).
Such kind of problems have been already considered by several authors ---
we refer the reader to a recent paper of Ustinov~\cite{Ustinov}, for example, where references to other papers
can be found. However we will not discuss these more general problems here.

\section{Proof of the theorem}

\subsection{Preparation} 

\indent

We assume that $q \le x^{\frac{2}{3}}$ because otherwise formula \eqref{30} is trivial.

\bigskip

The sum $S_{q, a}(x)$ defined by \eqref{10} is obviously equal to the number of pairs of integers
$u, v$ satisfying 
\begin{equation} \label{40}
  u^2 + v^2 \le x , \qquad u^2 + v^2 \equiv a \; (q) .
\end{equation}
Therefore
\begin{equation} \label{50}
  S_{q, a}(x) = 4 S' + 4 S^{\prime\prime} + O (1) ,
\end{equation}
where $S'$ is the number of pairs of positive integers $u, v$ satisfying \eqref{40}
and respectively $S^{\prime\prime} $ is the number of positive integers $u \le \sqrt{x} $ such that
\begin{equation} \label{60}
    u^2 \equiv a \; (q) .
\end{equation}
(From this point onwards by $u$ and $v$ we denote natural numbers only).

\bigskip

It is clear that
\begin{equation} \label{70}
  S^{\prime\prime} =  \frac{ \omega_a(q) }{q} \sqrt{x} + O \left(  \omega_a(q) \right)  ,
\end{equation}
where 
\begin{equation} \label{80}
  \omega_a(q) = \# \{ \;  1 \le \alpha \le q \; : \; \alpha^2 \equiv a \; (q)  \; \} .
\end{equation}
We note that this function satisfies
\begin{equation} \label{85}
  \omega_a(q) \ll (q, a)^{\frac{1}{2}} \tau(q) .
\end{equation}
This can be proved in a simple elementary way and we leave the verification to the reader.

\bigskip

The sum $S'$ can be written in the form
\begin{equation} \label{90}
  S' = 2 S_1 - S_2 .
\end{equation}
Here $S_1$ is the number of pairs of natural numbers $u, v$ satisfying the congruence
\begin{equation} \label{100}
 u^2 + v^2 \equiv a \; (q)
\end{equation}
and such that
\begin{equation} \label{110}
   u \le \sqrt{x/2} , \qquad  v \le \sqrt{x - u^2} .
\end{equation}
Respectively $S_2$ is the number pairs $u, v$ satisfying \eqref{100} and also
\begin{equation} \label{120}
   u  \le \sqrt{x/2} , \qquad      v  \le \sqrt{x/2} .
\end{equation}

\bigskip

We divide $S_1$ and $S_2$ into parts according to the congruence classes of $u$ and $v$ modulo $q$.
From this point we write for simplicity
\begin{equation} \label{125}
   \sum_{\alpha, \beta} \qquad \text{for} \qquad 
   \sum_{\substack{1 \le \alpha, \beta \le q \\ \alpha^2 + \beta^2 \equiv a \; (q)} } .
\end{equation}
In this notation we have
\begin{equation} \label{130}
  S_1 = \sum_{\alpha, \beta} G_1(\alpha, \beta) , \qquad S_2 = \sum_{\alpha, \beta} G_2(\alpha, \beta)
\end{equation}
where $ G_1(\alpha, \beta) $ is the number or pairs $u, v$ satisfying \eqref{110} and also
\begin{equation} \label{140}
  u \equiv \alpha \; (q) ,  \qquad v \equiv \beta \; (q) 
\end{equation}
and respectively $ G_2(\alpha, \beta) $ is the number of pairs $u, v$
with \eqref{120} and \eqref{140}.

\bigskip

Consider $G_1(\alpha, \beta)$. 
Denote as usual by $[y]$ and $\{y\}$ the integer part and the fractional part of $y$ and let
\begin{equation} \label{150}
  \rho(y) = \frac{1}{2} - \{ y \} .
\end{equation}
We use that for any $y \ge 0$ we have
\begin{equation} \label{155}
  \sum_{\substack{ u \le y  \\ u \equiv \gamma \; (q)} } 1
   = \left[ \frac{y - \gamma}{q} \right] - \left[ \frac{- \gamma}{q} \right] 
   = \frac{y}{q} + \rho \left( \frac{y - \gamma}{q} \right) - \rho \left( \frac{- \gamma}{q} \right) .
\end{equation}
Hence for any
$u$ satisfying the first condition in \eqref{110} the number of integers $v$ 
satisfying the second of these conditions and also the second congruence from \eqref{140} 
equals 
\[ 
  \frac{\sqrt{x - u^2}}{q} + \rho \left( \frac{\sqrt{x - u^2} - \beta}{q} \right) - 
  \rho \left( \frac{ - \beta}{q} \right) .
\]
Therefore
\[
 G_1(\alpha, \beta) = \frac{1}{q} \sum_{\substack{u \le \sqrt{x/2} \\ u \equiv \alpha \; (q)}} \sqrt{x - u^2}
 + \sum_{\substack{u \le \sqrt{x/2} \\ u \equiv \alpha \; (q)}} \rho \left( \frac{\sqrt{x - u^2} - \beta}{q} \right) 
 - \sum_{\substack{u \le \sqrt{x/2} \\ u \equiv \alpha \; (q)}} \rho \left( \frac{ -\beta}{q} \right) .
\]
We substitute this expression for $G_1$ in the first formula in \eqref{130} and we get
\begin{equation} \label{160}
  S_1 = \frac{1}{q} S_1^{(0)} + S_1^{(1)} - S_1^{(2)} ,
\end{equation}
where
\begin{align} 
  S_1^{(0)} 
  & = 
  \sum_{\alpha, \beta} \sum_{\substack{u \le \sqrt{x/2} \\ u \equiv \alpha \; (q)}} \sqrt{x - u^2} ,
\label{170}   \\
 & \notag \\
  S_1^{(1)} 
  & = 
  \sum_{\alpha, \beta}
   \sum_{\substack{u \le \sqrt{x/2} \\ u \equiv \alpha \; (q)}} \rho \left( \frac{\sqrt{x - u^2} - \beta}{q} \right) ,
\label{180} \\
& \notag \\
  S_1^{(2)} 
  & = 
   \sum_{\alpha, \beta} 
   \sum_{\substack{u \le \sqrt{x/2} \\ u \equiv \alpha \; (q)}} \rho \left( \frac{ - \beta}{q} \right) .
\label{190}
\end{align}
We proceed with $G_2$ in the same manner and using the second formula in \eqref{130} we get
\begin{equation} \label{200}
  S_2 = \frac{\sqrt{x/2}}{q} S_2^{(0)} + S_2^{(1)} - S_1^{(2)} ,
\end{equation}
where $S_1^{(2)}$ is specified by \eqref{190} and
\begin{align} 
  S_2^{(0)} 
  & = 
  \sum_{\alpha, \beta} \sum_{\substack{u \le \sqrt{x/2} \\ u \equiv \alpha \; (q)}} 1 ,
\label{210}   \\
 & \notag \\
  S_2^{(1)} 
  & = 
  \sum_{\alpha, \beta}
   \sum_{\substack{u \le \sqrt{x/2} \\ u \equiv \alpha \; (q)}} \rho \left( \frac{\sqrt{x/2} - \beta}{q} \right) .
\label{220} 
\end{align}
From \eqref{50}, \eqref{70}, \eqref{85}, \eqref{90}, \eqref{160} and \eqref{200} 
we obtain
\begin{align} 
    S_{q, a} (x) 
    & = 
    \frac{8}{q} S_1^{(0)} \, + \, 8 S_1^{(1)} \, - \, 4 S_1^{(2)} \,
  - \, 4 \frac{\sqrt{x/2}}{q} S_2^{(0)} 
    \notag \\
    & \notag \\
    & \qquad \qquad
        - \, 4 S_2^{(1)} \,
       + 4 \frac{\omega_a(q)}{q} \sqrt{x} \, + \, O \left( (q, a)^{\frac{1}{2}} \tau(q) \right) .
 \label{240}
\end{align}

\bigskip

We estimate the sum $S_1^{(1)}$ in sections \ref{S1_a} -- \ref{S1_d}. This is the most difficult part of the proof.
In section \ref{end} we evaluate $S_1^{(0)}$, $S_1^{(2)}$, $S_2^{(0)}$ and $S_2^{(1)}$.
Finally in section \ref{endoftheproof} we collect all results together and prove
the asymptotic formula \eqref{30}.

\subsection{Estimation of $S_1^{(1)}$ --- beginning.} \label{S1_a} 

\indent

It is well known that for any integer $M\ge 2$ 
the function $\rho (y)$ defined by \eqref{150} can be written in the form
\begin{equation} \label{250}
  \rho(y) = \sum_{1 \le |n| \le M } \frac{e(ny)}{2 \pi i n} + 
  O \left( \min \left( 1, \frac{1}{M ||y||} \right)   \right) ,
\end{equation}
where 
$||y||$ is the distance from $y$ to the nearest integer, 
$e(y) = e^{2 \pi i y}$, the constant in the $O$-symbol is absolute and
where, as usual, we assume that $\min \left( 1 , \frac{1}{0} \right) = 1$.
We also have
\begin{equation} \label{260}
  \min \left( 1, \frac{1}{M ||y||} \right)  = \sum_{n \in \mathbb Z} c_n e(ny) ,
\end{equation}
where
\begin{equation} \label{270}
  c_n \ll 
    \begin{cases}
      M^{-1} \log M  \quad & \text{for all}  \; n , \\
      M n^{-2}  & \text{for} \; n \not= 0 .
    \end{cases}
\end{equation}
For the proofs of \eqref{250} -- \eqref{270} we refer the reader to \cite[Chapter 2]{Hooley}.

\bigskip

We take an integer $M \ge 2$, which we shall choose later,
and using \eqref{180}, \eqref{250} we obtain
\begin{equation} \label{280}
  S_1^{(1)} = \sum_{\alpha, \beta}
   \sum_{\substack{u \le \sqrt{x/2} \\ u \equiv \alpha \; (q)}} \;
   \sum_{1 \le |n| \le M} 
   \frac{1}{2 \pi i n } \;
   e \left( \frac{\sqrt{x - u^2} - \beta}{q} \, n \right) 
   + O \left( \Delta \right) ,
\end{equation}
where
\begin{equation} \label{290}
  \Delta =  \sum_{\alpha, \beta} \sum_{\substack{u \le \sqrt{x/2} \\ u \equiv \alpha \; (q)}} 
  \min \left( 1 , M^{-1} \left| \left|   \frac{\sqrt{x - u^2} - \beta}{q} \right| \right|^{-1} \right) .
\end{equation}
Now we apply \eqref{260} to get
\begin{equation} \label{300}
  \Delta =  \sum_{\alpha, \beta} \; \sum_{\substack{u \le \sqrt{x/2} \\ u \equiv \alpha \; (q)}} \;
  \sum_{n \in \mathbb Z} c_n \; e \left( \frac{\sqrt{x - u^2} - \beta}{q} \, n \right) =
   \sum_{n \in \mathbb Z} c_n \; \mathcal F_n ,
\end{equation}
where
\begin{equation} \label{310}
 \mathcal F_n = \sum_{\alpha, \beta} \; \sum_{\substack{u \le \sqrt{x/2} \\ u \equiv \alpha \; (q)}} 
 e \left( \frac{\sqrt{x - u^2} - \beta}{q} \, n \right) .
\end{equation}
From \eqref{270}, \eqref{280}, \eqref{300} and \eqref{310} we find
\begin{equation} \label{320}
  S_1^{(1)} \ll \frac{\log M}{M} |\mathcal F_0| +
  \log M \, \sum_{1 \le |n| \le M}  |n|^{-1} \, |\mathcal F_n| 
  + M \, \sum_{|n| > M}  n^{-2} \, |\mathcal F_n| .
\end{equation}

\bigskip

Consider the sum $\mathcal F_n$. We use the elementary identity
\begin{equation} \label{325}
  \sum_{h \; (q)} e \left( \frac{hm}{q} \right) =
  \begin{cases}
    q \quad & \text{for} \quad q \mid m , \\
    0       & \text{otherwise} 
  \end{cases}
\end{equation}
(the summation is taken over all residue classes $h \pmod{q}$) 
and we find
\begin{align}
  \mathcal F_n 
  & = 
  \sum_{\alpha , \beta} \; \sum_{u \le \sqrt{x/2}} 
    e \left( \frac{\sqrt{x - u^2} - \beta}{q} n \right) \; \frac{1}{q} \sum_{h \; (q)} 
    e \left( \frac{h(u - \alpha)}{q} \right) 
    \notag \\
    & \notag \\
    & = 
    \frac{1}{q} \sum_{h \; (q)} \mathcal H_{h, n} \, \mathcal T_{h, n} ,
   \label{330}
\end{align}
where
\begin{align} 
   \mathcal H_{h, n} 
     & = 
     \mathcal  H_{h, n} (q, a) = \sum_{\alpha, \beta} e \left( \frac{-\alpha h - \beta n }{q} \right) ,
   \label{340} \\
   & \notag \\
     \mathcal T_{h, n} 
     & = 
     \mathcal  T_{h, n} (q) = \sum_{u \le \sqrt{x/2}} e (f(u)) , 
        \label{350} \\
     & \notag \\
        f(u) 
        & = \left(n \sqrt{x-u^2} +h  u  \right) \, q^{-1} .
        \label{355}
\end{align}
From \eqref{330} we obtain
\begin{equation} \label{360}
  |\mathcal F_n| \le \frac{1}{q} \sum_{|h| \le q/2 } |\mathcal H_{h, n} | \, | \mathcal T_{h, n} | ,
\end{equation}
so to proceed further we have to estimate the sums $\mathcal H_{h, n}$ and $\mathcal T_{h, n}$.

\subsection{Estimation of $\mathcal H_{h, n}$.} \label{S1_b}

\subsubsection{Preparation.} \label{beginning}

\indent

In section \ref{S1_b} we establish that
\begin{equation} \label{361}
  |\mathcal H_{h, n}(q, a)| \le 4 q^{\frac{1}{2}} \, \tau^2(q) \, (q, h, n)^{\frac{1}{2}} \, 
  (q, a, h^2 + n^2)^{\frac{1}{2}} ,
\end{equation}
where as usual $(u, v, w)$ is the greatest common divisor of $u, v, w$.
We note that in the proof of \eqref{30} we actually use only the following consequence of \eqref{361}:
\begin{equation} \label{430}
   | \mathcal H_{h, n} (q, a) | \le 4 q^{\frac{1}{2}} \, \tau^2(q) \, (q, h, n)^{\frac{1}{2}} \, 
   \left( q, a \right)^{\frac{1}{2}} .
\end{equation}

\bigskip

To establish \eqref{361} we first express $\mathcal H_{h, n}$ by means of the Gauss sum
\begin{equation} \label{362}
  S(q; k, m) = \sum_{\alpha \; (q)} e \left( \frac{k \alpha^2 + m \alpha}{q} \right) .
\end{equation}
Using \eqref{125}, \eqref{325} and \eqref{340} we get
\begin{align}
  \mathcal H_{h, n}(q, a) 
   & = 
   \sum_{\alpha, \beta }
  e \left( \frac{ \alpha h + \beta n}{q} \right) 
  = \sum_{\alpha \; (q) }  \sum_{\beta \; (q)}
  e \left( \frac{ \alpha h + \beta n}{q} \right)
  \frac{1}{q} \sum_{k \; (q)} e \left( \frac{k(\alpha^2 + \beta^2 - a)}{q} \right) 
   \notag \\
   & \notag \\
   & = 
   \frac{1}{q} \sum_{k \; (q)} 
     e \left( \frac{- a k}{q} \right)   \; S(q; k, h) \; S(q; k, n) .
     \label{390}
\end{align}
We know that 
\begin{equation} \label{400}
  S(q; k, n) = 
  \begin{cases} 
      d  S \left( q/d; k/d, n/d \right) \quad & \text{for} 
           \quad d \mid n , \qquad \text{where} \quad d = (k, q) , \\
      0  & \text{for} \quad d \nmid n .
  \end{cases}
\end{equation}
A proof of this relation can be found in \cite[Section 6]{Esterman}.
Therefore from \eqref{390} we find
\begin{align}
  \mathcal H_{h, n} (q, a) 
    & = 
    \frac{1}{q} \sum_{ d \mid q} 
  \sum_{\substack{k \; (q) \\ (k, q) = \frac{q}{d}}} e \left( \frac{- a k}{q} \right) \; S(q; k, h) \; S(q; k, n) 
   \notag \\
   & \notag \\
   & =
    q \sum_{\substack{ d \mid q \\ \frac{q}{d} \mid (h, n)} } d^{-2} \;
  B(d; a, hd/q, nd/q) ,
   \label{410}
\end{align}
where
\begin{equation} \label{412}
   B(d; a, m, t) =
   {\sum_{l \; (d) }}^* e \left( \frac{- a l}{d} \right) 
   \; S(d; l, m )    \;  S(d; l, t ) 
\end{equation}
(the asterisk means that the summation is restricted to a reduced system of residues).

\bigskip

Next we use other properties of the Gauss sum to represent 
the quantity $B$ by means of the Kloosterman sum
\begin{equation} \label{370}
  K(q; k, n) = {\sum_{\alpha \; (q)}}^* e \left( \frac{k \alpha + n \overline{\alpha}}{q} \right) ,
\end{equation}
where $\overline{\alpha}$ is the inverse of $\alpha$ modulo $q$,
and then we apply A.Weil's bound
\begin{equation} \label{380}
  |K(q; k, n)| \le q^{\frac{1}{2}} \, \tau(q) \, (q, k, n)^{\frac{1}{2}} .
\end{equation}
A proof of \eqref{380} is available in \cite[Chapter~11]{IwKo}.

\bigskip

The calculations are simpler if $q$ is odd but certain technical complications arise in the general case.
That is why we consider separately the cases $2 \nmid q$ and $q = 2^{\theta}$ and establish 
corresponding versions of \eqref{361}.
Finally we note that $\mathcal H_{h, n}(q, a)$ has a multiplicative property with respect to $q$ and 
prove \eqref{361}.

\subsubsection{The case $2 \nmid q$.} \label{2nmq}

\indent

In this section we prove that
\begin{equation} \label{420}
   | \mathcal H_{h, n} (q, a) | \le q^{\frac{1}{2}} \, \tau^2(q) \, (q, h, n)^{\frac{1}{2}} \, 
   \left( q, a, h^2 + n^2 \right)^{\frac{1}{2}} \qquad \quad \text{for} \qquad \quad 2 \nmid q .
\end{equation}
Having in mind \eqref{410} we see that we have to estimate
the quantity $B$ defined by \eqref{412}
for any $d \mid q$.

\bigskip

We first note that the Gauss sum satisfies
\begin{equation} \label{415}
  S(q; k, m) = 
  e \left( \frac{- \overline{(4k)} \, m^2}{q} \right) \; \left( \frac{k}{q} \right) \; S(q; 1) 
  \qquad \text{for} \qquad (q, 2k)=1 ,
\end{equation}
where $\left( \frac{k}{q} \right)$ is the Jacobi symbol and $S(q; 1) = S(q; 1, 0)$.
A proof of \eqref{415} can be found in \cite[Section 6]{Esterman}.
Using \eqref{412}, \eqref{370}, \eqref{415} and bearing in mind that
$ 2 \nmid d $ for any $d \mid q$ we get
\[
  B(d; a, m, t) = S^2(d; 1) \, K (d; a, \overline{4} (m^2 + t^2))  ,
\]
where $\overline 4$ is the inverse of $4$ modulo $d$.
It is well known that
\[
  |S(d; 1)|^2 = d \qquad \text{for} \qquad 2 \nmid d
\]
(a proof is available for example \cite[Chapter 7]{Hua}). From this formula and \eqref{380} we find
\[
 |B(d; a, m, t)| \le d^{\frac{3}{2}} \, \tau(d) \, (d, a, m^2 + t^2)^{\frac{1}{2}} .
\]
We substitute this bound for $B$ in \eqref{410} 
(with $m = hd/q$ and $t = nd/q$)
to get
\[
  |\mathcal H_{h, n}(q, a)| \le  
   q \, \tau(q) \, (q, a, h^2 + n^2)^{\frac{1}{2}} \, \mathfrak X , \qquad 
   \mathfrak X =
   \sum_{\substack{d \mid q \\ \frac{q}{d} \mid (h, n)} } d^{- \frac{1}{2}} .
\]
It is obvious that
\[
  \mathfrak X = \sum_{\delta \mid (q, h, n)} \left( \frac{q}{\delta} \right)^{-\frac{1}{2}} \le 
  q^{-\frac{1}{2}} \, \tau(q) \, (q, h, n)^{\frac{1}{2}} ,
\]
hence we obtain \eqref{420}.

\subsubsection{The case $q = 2^{\theta}$.}

In this section we prove that
\begin{equation} \label{380.6}
  \left| \mathcal H_{h, n} \left(2^{\theta}, a \right) \right| 
  \le 4 \left( 2^{\theta} \right)^{\frac{1}{2}} \, \tau^2 \left( 2^{\theta} \right) \, 
   \left( 2^{\theta}, n, h \right)^{\frac{1}{2}} \, 
   \left( 2^{\theta}, a, h^2 + n^2 \right)^{\frac{1}{2}}  .
\end{equation}
First we establish that the sum $B$ given by \eqref{412} satisfies
\begin{equation} \label{382}
   \left| B \left(2^{\nu}, a, m, t \right) \right| 
  \le 4 \left( 2^{\nu} \right)^{\frac{3}{2}} \, \tau \left( 2^{\nu} \right) \, 
   \left( 2^{\nu}, a, m^2 + t^2 \right)^{\frac{1}{2}}  .
\end{equation}

\bigskip

The inequality \eqref{382} is obvious for $\nu < 2$, so we may assume that $\nu \ge 2$.
Under this assumption and if $2 \nmid r$ the Gauss sum satisfies
\begin{equation} \label{383}
  S \left( 2^{\nu}, r, t \right) =
  \begin{cases}
     e \left( \frac{- \overline{r} \; (t/2)^2 }{2^{\nu}} \right) \, 2^{\frac{\nu+1}{2}} \;
     \frac{1 + i^r}{\sqrt{2}} , 
     \qquad & \text{for} \quad 2 \mid t , \quad 2 \mid \nu  ,  \\
      e \left( \frac{- \overline{r} \; (t/2)^2 }{2^{\nu}} \right) \, 2^{\frac{\nu+1}{2}} \;
       e \left( \frac{r}{8} \right) , 
     \qquad & \text{for} \quad 2 \mid t , \quad 2 \nmid \nu  ,  \\
     0 & \text{for} \quad 2 \nmid t .
  \end{cases}
\end{equation}
For a proof we refer the reader to \cite[Section 6]{Esterman} and \cite[Chapter 7]{Hua}.

\bigskip

From \eqref{412} and \eqref{383} it follows that \eqref{382} holds if
$2 \nmid m$ or $2 \nmid t$.

\bigskip

Consider now the case when $2 \mid m $ and $2 \mid t$. 
If we have also $2 \mid \nu$ then \eqref{412} and \eqref{383} imply
\[
  B \left( 2^{\nu}, a, m, t \right) = i 2^{\nu + 1} \, { \sum_{l \; (2^{\nu} )} }^*
  (-1)^{\frac{l-1}{2}} \, 
  e \left( \frac{ -l a - \overline{l} 
  \left( \left(\frac{m}{2} \right)^2 + \left( \frac{t}{2} \right)^2 \right) }{2^{\nu}} \right) 
\]
and having in mind \eqref{370} we get
\begin{equation} \label{384}
   B \left( 2^{\nu}, a, m, t \right) = 
   2^{\nu + 1} \, K \left( 2^{\nu}; 2^{\nu - 2} - a , - \frac{m^2 + t^2}{4} \right)  .
\end{equation}
Now we apply \eqref{380} and obtain
\[
   | B \left( 2^{\nu}, a, m, t \right) | \le 2 \, \left( 2^{\nu} \right)^{\frac{3}{2}} \, 
   \tau \left( 2^{\nu} \right) \, \left( 2^{\nu} , 2^{\nu - 2} - a , \frac{ m^2 + t^2 }{4} \right)^{\frac{1}{2}} ,
\]
which implies \eqref{382}. 

\bigskip

If $2 \nmid \nu$ then using \eqref{412}, \eqref{370} and \eqref{383} we see that 
\eqref{384} is true again, hence \eqref{382} holds also in this case.

\bigskip

Now we proceed as in section \ref{2nmq} and using \eqref{382} we obtain \eqref{380.6}.

\subsubsection{The estimate for $\mathcal H_{h, n}$.}

\indent

We are now in a position to prove \eqref{361} in the general case. We have the following identity
\begin{equation} \label{381}
  \mathcal H_{h, n}(q_1 q_2, a) = \mathcal H_{h\overline{q_2}, \, n\overline{q_2}}(q_1, a) \;
    \mathcal H_{h\overline{q_1}, \, n \overline{q_1}}(q_2, a) ,
    \qquad \text{for} \qquad (q_1, q_2)=1 ,
\end{equation}
where $\overline{q_2}$ and $\overline{q_1}$ denote the inverses of $q_2$ and $q_1$ 
modulo $q_1$ and respectively $q_2$. The proof of \eqref{381} is standard and we leave it to
the reader.
Now we represent $q = q_1 q_2$, where $2 \nmid q_1$ and $q_2 = 2^{\theta}$, then we apply 
\eqref{420}, \eqref{380.6} and \eqref{381} and obtain \eqref{361}.

\subsection{Estimation of $\mathcal T_{h, n}$ and $\mathcal F_n$.} \label{S1_c}

\subsubsection{The sum $\mathcal T_{h, n}$.}

\indent

In this section we assume that $|h| \le \frac{q}{2}$ and establish that
\begin{equation} \label{590}
\mathcal T_{h, n}  \ll 
    \min( \sqrt{x}, q |h|^{-1} ) \qquad \text{if} \qquad n = 0 \qquad \text{or} \qquad 0 < 2|n| \le |h|
\end{equation}
and
\begin{equation} \label{595}
  \mathcal T_{h, n}  \ll 
      x^{\frac{1}{4}} \left( |n|^{\frac{1}{2}} q^{-\frac{1}{2}} + |n|^{-\frac{1}{2}} q^{\frac{1}{2}} \right) 
      \qquad \text{for} \qquad n \not= 0  .
\end{equation}

\bigskip

If $n=0$ then the inequality \eqref{590} follows from the well-known estimate for
linear exponentials sums (see \cite[Chapter 2.1]{GraKol}).

\bigskip

Suppose now that $n \not= 0$.
The function $f(u)$ defined by \eqref{355} satisfies
\[
  f'(u)  = 
  \left( - n u (x - u^2)^{- \frac{1}{2}} + h \right) q^{-1} , 
  \qquad f^{\prime\prime}(u)  = 
  - n  x (x-u^2)^{-\frac{3}{2}} \, q^{-1} .
\]
In the case $ 0 < 2|n| \le |h| \le \frac{q}{2}  $ we use the first formula
and find
\[
  \frac{|h|}{2 q} \le |f'(u)| \le \frac{3|h|}{2q} \le \frac{3}{4}
   \qquad \text{for} \qquad 0 \le u \le \sqrt{x/2} .
\]
Then we apply \cite[Theorem 2.1]{GraKol} and we obtain \eqref{590}.

\bigskip

Further, it is clear that
\[
  |f^{\prime\prime}(u)| \asymp |n| \, q^{-1} \, x^{-\frac{1}{2}} \qquad \text{for} \qquad 0 \le u \le \sqrt{x/2} .
\]
Now we apply \cite[Theorem 2.2]{GraKol} and find
\[
  \mathcal T_{n, h} \ll x^{\frac{1}{2}} \left( |n| \, q^{-1} \, x^{-\frac{1}{2}} \right)^{\frac{1}{2}}
    + \left( |n| \, q^{-1} \, x^{-\frac{1}{2}} \right)^{-\frac{1}{2}} ,
\]
which gives \eqref{595}.

\subsubsection{The sum $\mathcal F_n$.}

\indent

If $n=0$ then we apply \eqref{360}, \eqref{430} and \eqref{590} to get 
\begin{align}
    \mathcal F_0 
    & \ll 
    q^{- \frac{1}{2}} \, \tau^2(q) \, (q, a)^{\frac{1}{2}} \,
    \sum_{|h| \le \frac{q}{2}} (q, h)^{\frac{1}{2}} \, \min (\sqrt{x} , q |h|^{-1})
     \notag \\
     & \notag \\
     & \ll
      \tau^2(q) \, (q, a)^{\frac{1}{2}} 
      \left( x^{\frac{1}{2}} + q^{\frac{1}{2} } \sum_{1 \le h \le q } \frac{(q, h)^{\frac{1}{2}}}{h} \right) .
      \notag
\end{align} 
However for any $y \ge 2$ we have
\begin{equation} \label{610}
   \sum_{ n \le y } \frac{(q, n)^{\frac{1}{2}}}{n} \le
   \sum_{\delta \mid q} \delta^{\frac{1}{2}} \, \sum_{\substack{n \le y \\ n \equiv 0 \; (\delta)} } \frac{1}{n}
   \le \sum_{\delta \mid q} \delta^{-\frac{1}{2}} \, \sum_{n \le \frac{y}{\delta} } \frac{1}{n}
   \ll \tau(q) \log y ,
\end{equation}
hence we get
\begin{equation} \label{620}
   \mathcal F_0 \ll x^{\frac{1}{2}} \,
   \tau^3(q) \, (q, a)^{\frac{1}{2}} \, \log x .
\end{equation}

\bigskip

In the case 
$n \not= 0$ we apply \eqref{360} and \eqref{430} to find
\[
\mathcal F_n  \ll q^{-\frac{1}{2}} \, \tau^2(q) \, (q, a)^{\frac{1}{2}} \, (q, n)^{\frac{1}{2}} \; 
  \mathfrak T_n , \qquad 
  \mathfrak T_n =  \sum_{|h| \le \frac{q}{2}} |\mathcal T_{h, n}| .
\]
However from \eqref{590} and \eqref{595} it follows that
\[
  \mathfrak T_n \ll
  \sum_{2|n| \le |h| \le \frac{q}{2}} q |h|^{-1} + \sum_{|h| \le \min \left( 2|n|, \frac{q}{2} \right)} 
  x^{\frac{1}{4}} \left( |n|^{\frac{1}{2}} q^{-\frac{1}{2}} + |n|^{-\frac{1}{2}} q^{\frac{1}{2}} \right)
  \ll q \log x +  x^{\frac{1}{4}} \, |n|^{\frac{1}{2}} \, q^{\frac{1}{2}} 
\]
and we get
\begin{equation} \label{600}
      \mathcal F_n \ll 
      \left( q^{\frac{1}{2}} + x^{\frac{1}{4}} |n|^{\frac{1}{2}} \right) \, 
      \tau^2(q) \, (q, a)^{\frac{1}{2}}\, (q, n)^{\frac{1}{2}} \, \log x  
      \qquad \quad \text{for} \qquad \quad n \not= 0 .
\end{equation}

\subsection{Estimation of $S_1^{(1)}$ --- end.} \label{S1_d}

\indent

We use \eqref{320}, \eqref{620} and \eqref{600} to find
\begin{align}
  S_1^{(1)} 
  & \ll 
     \left( 
   x^{\frac{1}{2}} M^{-1} 
    + 
   \left( q^{\frac{1}{2}} + x^{\frac{1}{4}} M^{\frac{1}{2}} \right) 
  \sum_{1 \le n \le M} \frac{(q, n)^{\frac{1}{2}}}{n}  
  + 
     M \sum_{n > M} \frac{(q, n)^{\frac{1}{2}}}{n^2} \, 
      \left( q^{\frac{1}{2}} + x^{\frac{1}{4}} n^{\frac{1}{2}} \right)  \right)    .   
      \notag \\
      & \notag \\
      & \qquad \qquad
         \times \, \tau^3(q) \, (q, a)^{\frac{1}{2}} \, \log x  \, \log M .
   \label{625} 
\end{align}       
However for any $y \ge 1 $ we have
\begin{equation} \label{630}
   \sum_{ n > y } \frac{(q, n)^{\frac{1}{2}}}{n^2} \le
   \sum_{\delta \mid q} \delta^{\frac{1}{2}} \, \sum_{\substack{n > y \\ n \equiv 0 \; (\delta)} } \frac{1}{n^2}
   \le \sum_{\substack{\delta \mid q \\ \delta \le y }}\delta^{-\frac{3}{2}} \, \sum_{n > \frac{y}{\delta} }   
     \frac{1}{n^2} + 
     \sum_{\substack{\delta \mid y \\ \delta > y }} \delta^{-\frac{3}{2}}
   \ll y^{-1} \,\tau(q)  
\end{equation}
and similarly
\begin{equation} \label{640}
  \sum_{n > y} \frac{(q, n)^{\frac{1}{2}}}{n^{\frac{3}{2}}}
  \ll y^{-\frac{1}{2}} \, \tau(q) .
\end{equation}
From \eqref{610}, \eqref{625} -- \eqref{640} it follows that
\[
  S_1^{(1)} \ll 
    \left( x^{\frac{1}{2}} M^{-1} + q^{\frac{1}{2}} + x^{\frac{1}{4}} M^{\frac{1}{2}} \right) 
    \, \tau^4(q) \, (q, a)^{\frac{1}{2}} \, \log x \, \log^2 M .
\]
We choose   $ M = \big[ x^{\frac{1}{6}} \, \big] $ and obtain
\begin{equation} \label{650}
  S_1^{(1)} \ll \left( q^{\frac{1}{2}} + x^{\frac{1}{3}}\right) 
  \,  \tau^4(q) \, (q, a)^{\frac{1}{2}} \, \log^4 x .
\end{equation}

\subsection{Evaluation of the sums $S_1^{(0)}$, $S_1^{(2)}$ $S_2^{(0)}$, $S_2^{(1)}$.} \label{end}

\subsubsection{The sum $S_1^{(0)}$.}

\indent

We use \eqref{80} and \eqref{125} to write the sum $S_1^{(0)}$ defined by \eqref{170} in the form
\[
  S_1^{(0)} = \sum_{u \le \sqrt{x/2}}  b_u \; \sqrt{x - u^2} , \qquad 
  b_u = \omega_{a - u^2}(q)
\]
and then apply Abel's transformation to get
\begin{equation} \label{670}
  S_1^{(0)} = \sqrt{x/2} 
  \sum_{u \le \sqrt{x/2}} b_u - \int_0^{\sqrt{x/2}}  \left( \sum_{u \le t} b_u \right) 
  \; \frac{d}{dt} \sqrt{x - t^2} \; d t .
\end{equation}
According to \eqref{125} and \eqref{155} we have
\[
   \sum_{u \le t} b_u  = \sum_{\alpha, \beta} \sum_{\substack{u \le t \\ u \equiv \alpha \; (q)}} 1
   = \sum_{\alpha, \beta} 
   \left( \frac{t}{q} + \rho \left( \frac{t - \alpha}{q} \right) - \rho \left( \frac{- \alpha}{q} \right)  \right) ,
\]
hence after certain simple calculations which we leave to the reader
and using \eqref{17}, \eqref{125} and \eqref{670} we find that
\begin{equation} \label{690}
  S_1^{(0)} 
     = \left( \frac{\pi}{8} + \frac{1}{4} \right) \, \frac{\eta_a(q)}{q} \, x
    + \sqrt{x/2} \; \mathfrak N    \; 
    - \sqrt{x} \; \mathfrak N_0 
      \; + \; \sum_{\alpha, \beta}  \Gamma_{\alpha} ,
\end{equation}            
where
\begin{equation} \label{695}
  \mathfrak N = 
    \sum_{\alpha, \beta} 
  \rho \left( \frac{\sqrt{x/2} - \alpha}{q} \right) ,  
  \qquad \mathfrak N_0 = \sum_{\alpha, \beta} \rho \left( \frac{- \alpha}{q} \right) 
\end{equation}
and
\begin{equation} \label{700}
     \Gamma_{\alpha} 
    = \int_{0}^{\sqrt{x/2}} \rho \left( \frac{t - \alpha}{q} \right) \, 
      \frac{t \; d t}{\sqrt{x - t^2}} .
\end{equation}

\bigskip

First we prove that
\begin{equation} \label{750}
 \mathfrak N_0  = \frac{1}{2} \, \omega_a(q) .
\end{equation}
Using \eqref{125} and \eqref{695} we can write
\begin{equation} \label{760}
   \mathfrak N_0  
   = \sum_{\beta \; (q)} \sum_{\substack{ 1 \le \alpha \le q \\ \alpha^2 \equiv a - \beta^2 \; (q) }} 
  \rho \left( \frac{-\alpha}{q} \right)
  =  \sum_{\beta \; (q)} \mathcal Y_{\beta} ,
\end{equation}
say. 

\bigskip

Consider the sum $\mathcal Y_{\beta}$. 
If $\beta^2 \not\equiv a \; (q)$ then we have $\mathcal Y_{\beta} = 0$. Indeed, in this case there is no term 
corresponding to $\alpha = q$; the term corresponding to $\alpha = q/2$ (if such exists) is equal to zero
and the other terms can be divided into couples 
$\rho \left( \frac{-\alpha}{q} \right) + \rho \left( \frac{\alpha - q}{q} \right) $, where
$1 \le \alpha < q/2$ and the sum of the terms of each such couple equals zero.

\bigskip

If $\beta^2 \equiv a \; (q)$ then we have $\mathcal Y_{\beta} = \frac{1}{2}$. 
Indeed, arguing as above we see that the contribution to $\mathcal Y_{\beta}$
from the terms corresponding to $1 \le \alpha < q$ vanishes. 
In the present case however there is a term corresponding to $\alpha = q$ and its contribution equals $\frac{1}{2}$.
This proves \eqref{750}.

\bigskip

We note that the above arguments imply also 
\begin{equation} \label{765}
    \left| \sum_{\alpha, \beta} \rho \left( \frac{-\alpha}{q} \right) \, \xi_{\beta} \right|
    \le 
    \omega_{a}(q) 
    \qquad \quad \text{for} \qquad \quad
    |\xi_{\beta}| \le 1 
\end{equation}
and
\begin{equation} \label{766}
     \sum_{\alpha, \beta} 
     \sin \left( \frac{2 \pi n \alpha}{q} \right) = 0 .     
\end{equation}

\bigskip

Consider now the integral $\Gamma_{\alpha}$.
The function $\rho(y)$ defined by \eqref{150} satisfies
\[
   \rho(y) = \sum_{n=1}^{\infty} \frac{\sin( 2 \pi n y)}{\pi n} \qquad \text{for} \qquad y \not\in \mathbb Z .
\]   
We insert this expression in \eqref{700} and change the order of summation and integration
(here we appeal to the dominated convergence theorem) and we find
\[
   \Gamma_{\alpha} = \sum_{n=1}^{\infty} \frac{1}{ \pi n} \, \int_{0}^{\sqrt{x/2}} 
   \sin \left( 2 \pi n \, \frac{t - \alpha}{q}  \right) 
   \; 
      \frac{t \; d t}{\sqrt{x - t^2}} .
\]
Hence using \eqref{766} we get
\begin{equation} \label{710}
   \sum_{\alpha, \beta} \Gamma_{\alpha} =  \sum_{n=1}^{\infty} \frac{1}{ \pi n}
   \, \mathcal D_n \, \mathcal E_n , 
\end{equation}
where
\[
  \mathcal D_n = \sum_{\alpha, \beta} \cos \left( \frac{2 \pi n \alpha}{q} \right) ,
  \qquad
  \mathcal E_n = \int_0^{\sqrt{x/2}}
  \sin \left( 2 \pi n \frac{t}{q} \right) 
  \;       \frac{t \; d t}{\sqrt{x - t^2}} .
\]
From \eqref{340} and \eqref{430} we find
\begin{equation} \label{720}
   \mathcal D_n \ll q^{\frac{1}{2}} \, \tau^2(q) \, (q, n)^{\frac{1}{2} } \, (q, a)^{\frac{1}{2}}
\end{equation}
and integrating by parts we find that
\begin{equation} \label{730}
 \mathcal E_n \ll \frac{q}{n} .
\end{equation}
Hence applying \eqref{630} (with $y=1$) and using 
\eqref{710} -- \eqref{730}
we obtain
\begin{equation} \label{740}
 \sum_{\alpha, \beta} \Gamma_{\alpha} \ll q^{\frac{3}{2}} \, \tau^3(q)  \, (q, a)^{\frac{1}{2}} .
\end{equation}

\bigskip

From \eqref{690}, \eqref{750} and \eqref{740} we obtain
\begin{equation} \label{770}
   S_1^{(0)} 
      = \left( \frac{\pi}{8} + \frac{1}{4} \right) \, \frac{\eta_a(q)}{q} \, x
    + \sqrt{x/2} \; \mathfrak N 
    - \frac{\sqrt{x}}{2} \omega_a(q)
     + O \left( q^{\frac{3}{2} } \, \tau^3(q) \, (q, a)^{\frac{1}{2}}  \right) .
\end{equation}
We note that there is no need to study the sum
$\mathfrak N $ because in the final expression for $S_{q, a}(x)$ the terms including it 
cancel each other.

\subsubsection{The sum $S_1^{(2)}$.}

\indent

From \eqref{85}, \eqref{155}, \eqref{190}, \eqref{750} and \eqref{765} we easily find
\begin{align}
 S_1^{(2)} 
   & = \sum_{\alpha, \beta} \rho \left( \frac{ - \beta}{q} \right) 
 \left( \frac{\sqrt{x/2}}{q} + 
 \rho \left(  \frac{\sqrt{x/2} - \alpha}{q} \right)  - \rho \left( \frac{- \alpha}{q} \right) \right) 
 \notag \\
 & \notag \\
   & =
    \frac{\omega_a(q)}{2 q} \sqrt{x/2}  \, +  \, O \left( (q, a)^{\frac{1}{2}} \tau(q) \right) .
   \label{780}
\end{align}

\subsubsection{The sum $S_2^{(0)}$.}

\indent

Using \eqref{17}, \eqref{85}, \eqref{125}, \eqref{155}, \eqref{210}, \eqref{695} and \eqref{750} we find
\begin{align}
   S_2^{(0)} 
    & = 
    \sum_{\alpha, \beta} \left( \frac{\sqrt{x/2}}{q} 
     + \rho  \left( \frac{\sqrt{x/2} - \alpha}{q} \right) - 
       \rho  \left( \frac{ - \alpha}{q} \right) 
     \right)
     \notag \\
     & \notag \\
     & =
     \frac{\eta_a(q)}{q} \sqrt{x/2} \, + \, \mathfrak N \, 
     - \frac{1}{2} \omega_a(q) .
   \label{790}
\end{align}

\subsubsection{The sum $S_2^{(1)}$.}

\indent

From \eqref{85}, \eqref{155}, \eqref{220}, \eqref{695} and \eqref{765} we get
\begin{align}
  S_2^{(1)} 
   & = 
   \sum_{\alpha, \beta} \rho \left( \frac{\sqrt{x/2} - \beta}{q} \right)
  \left( \frac{\sqrt{x/2}}{q} + 
 \rho \left(  \frac{\sqrt{x/2} - \alpha}{q} \right)  - \rho \left( \frac{- \alpha}{q} \right) \right) 
   \notag \\
   & \notag \\
   & =
   \frac{\sqrt{x/2}}{q} \, \mathfrak N  \, 
   + \, \mathfrak D  \, + O \left( (q, a)^{\frac{1}{2}} \, \tau(q) \right) ,
   \label{800}
\end{align}
where
\begin{equation} \label{810}
   \mathfrak D = \sum_{\alpha, \beta} 
    \rho \left(  \frac{\sqrt{x/2} - \alpha}{q} \right) 
     \rho \left(  \frac{\sqrt{x/2} - \beta}{q} \right) .
\end{equation}

\bigskip

Consider the sum $\mathfrak D$. We take an integer $M_1 \ge 2$ and apply \eqref{250} with $M=M_1$ to get
\[
  \mathfrak D   =
    \sum_{\alpha, \beta} \rho \left(  \frac{\sqrt{x/2} - \alpha}{q} \right) 
     \sum_{1 \le |n| \le M_1} \frac{1}{2 \pi i n} \,
      e \left(  \frac{\sqrt{x/2} - \beta}{q} \, n \right) 
      + O \left( \Delta_1 \right) ,
\]
where
\begin{equation} \label{820}
 \Delta_1 = \sum_{\alpha, \beta}  \min \left( 1, M_1^{-1} 
 \left|\left| \frac{\sqrt{x/2} - \beta }{q} \right|\right|^{-1} \right) .
\end{equation}
Applying \eqref{250} again and having in mind \eqref{340} and \eqref{820} we find
\begin{align}
   \mathfrak D   
   & =
    \sum_{\alpha, \beta} \;
     \sum_{1 \le |m| \le M_1} \frac{1}{2 \pi i m} \,
      e \left(  \frac{\sqrt{x/2} - \alpha}{q} \, m \right) 
       \sum_{1 \le |n| \le M_1} \frac{1}{2 \pi i n} \,
      e \left(  \frac{\sqrt{x/2} - \beta}{q} \, n \right)  
      \notag \\
      & \notag \\
      & 
      \qquad \qquad \qquad \qquad \qquad \qquad + 
      O \left( \Delta_1 \log M_1 \right) 
      \notag \\
      & \notag \\
      & =
      \sum_{1 \le |m|, |n| \le M_1} \frac{e \left( (m+n) \, q^{-1} \sqrt{x/2}  \right)}{(2 \pi i)^2 mn} \;
      \mathcal H_{m, n}  \; + \; O \left( \Delta_1 \log M_1 \right) .
      \label{830}
\end{align}
Next we use 
\eqref{260}  (with $M=M_1$) to write
the sum $\Delta_1$ defined by \eqref{820} in the form
\begin{equation} \label{840}
  \Delta_1 = \sum_{\alpha, \beta} \; \sum_{n \in \mathbb Z} c_n \; e \left( \frac{\sqrt{x/2} - \beta}{q} \, n \right)
  = \sum_{n \in \mathbb Z} c_n \, e \left( \frac{\sqrt{x/2}}{q} \, n \right) \, \mathcal H_{0, n} .
\end{equation}
From \eqref{17}, \eqref{270}, \eqref{430}, \eqref{610}, \eqref{630}, \eqref{830} and \eqref{840} we find
\begin{align}
  \mathfrak D 
    & \ll 
   \sum_{1 \le |m|, |n| \le M_1} \frac{|\mathcal H_{m, n}|}{|mn|} + 
   \log^2 M_1
   \left( 
   \frac{\eta_a (q) }{M_1} 
   + 
    \sum_{1 \le |n| \le M_1} \frac{|\mathcal H_{0, n}| }{|n|} 
   + M_1 \sum_{|n| > M_1} \frac{|\mathcal H_{0, n}| }{n^2}  \right) .
   \notag \\
   & \notag \\
     & \ll
    q^{\frac{1}{2}} \, \tau^2 (q) \, (q, a)^{\frac{1}{2}} \, 
    \left(
   \sum_{ 1 \le m, n \le M_1} \frac{(q, m, n)^{\frac{1}{2}}}{mn} +
      \sum_{1 \le n \le M_1} \frac{(q, n)^{\frac{1}{2}}}{n} + M_1 \sum_{n > M_1} \frac{(q, n)^{\frac{1}{2}}}{n^2} \right)
      \, \log^2 M_1 
      \notag \\
      & \notag \\
   & \qquad \qquad \qquad
         +        \frac{\log^2 M_1  }{M_1} \eta_a (q) 
           \notag \\
      & \notag \\
   & \ll
     q^{\frac{1}{2}} \, \tau^3 (q) \, (q, a)^{\frac{1}{2}} \, \log^4 M_1 
      +        \frac{\log^2 M_1  }{M_1} \eta_a (q) .
   \notag
\end{align}
Now we choose $M_1 = q^2 $ and having in mind 
\eqref{17} we get
\begin{equation} \label{850}
 \mathfrak D \ll q^{\frac{1}{2}} \, \tau^3 (q) \, (q, a)^{\frac{1}{2}} \log^4 x .
\end{equation}
From \eqref{800} and \eqref{850} we obtain
\begin{equation} \label{860}
 S_2^{(1)} = \frac{\sqrt{x/2}}{q} \, \mathfrak N \, + \, 
   O \left( q^{\frac{1}{2}} \, \tau^3 (q) \, (q, a)^{\frac{1}{2}} \log^4 x  \right)
\end{equation}

\subsection{The end of the proof.} \label{endoftheproof}

\indent

It remains to collect together \eqref{240}, \eqref{650}, \eqref{770}, \eqref{780}, \eqref{790}
and \eqref{860} and we establish \eqref{30}, which proves the theorem.

\bigskip
\bigskip

\vbox{
\hbox{Faculty of Mathematics and Informatics}
\hbox{Sofia University ``St. Kl. Ohridsky''}
\hbox{5 J.Bourchier, 1164 Sofia, Bulgaria}
\hbox{ }
\hbox{Email: dtolev@fmi.uni-sofia.bg}}

\end{document}